\documentclass[12pt]{amsart} 
\usepackage{verbatim, amssymb}

\def\AAP{Ann.\ Appl.\ Prob.\ }
\def\BER{Bernoulli\ }%
\def\CMMP{USSR Comput.\ Math.\ Math.\ Phys.\ }
\def\JC{J.\ Complexity\ }
\def\JFAA{J.\ Fourier Anal.\ Appl.\ }%
\def\JMP{J.\ Math.\ Phys.\ }
\def\JTP{J.\ Theoret.\ Probab.\ }
\def\LNM{Lect.\ Notes in Math.\ }
\def\PTRF{Probab.\ Theory Relat.\ Fields\ }
\def\SPL{Stat.\ Prob.\ Letters\ }
\def\TMM{Transl.\ Math.\ Monogr.\ }%

\newcommand{\R}{{\mathbb R}}  
\newcommand{\N}{{\mathbb N}} 
\newcommand{\eps}{\varepsilon}
\newcommand{\dist}{\operatorname{dist}}

\newcommand{\XX}{{\mathfrak X}}
\newcommand{\YY}{{\mathfrak Y}}
\newcommand{\SSS}{{\mathbb S}} 
\newcommand{\deter}{\mathrm{det}} 
\newcommand{\ran}{\mathrm{ran}} 
 
\newcommand{\cost}{\operatorname{cost}} 
\newcommand{\card}{\operatorname{card}} 
\newcommand{\Sh}{\widehat{S}} 
\newcommand{\Xh}{\widehat{X}}
 
\newcommand{\Wt}{\widetilde{W}} 
\newcommand{\Xb}{\overline{X}} 
\newcommand{\A}{{\mathfrak A}}
\newcommand{\qr}{q^{(r)}} 
\newcommand{\qe}{q^{(1)}}                 
\newcommand{\de}{d^{(1)}}               
\newcommand{\grad}{\nabla}
\newcommand{\nice}[2]{\phantom{#2} #1, #2}

\renewcommand{\P}{{\mathbb P}}  
\newcommand{\E}{{\mathbb E}}
\newcommand{\V}{\mathbb{V}} 

\theoremstyle{plain}
\newtheorem{theorem}{Theorem}
\newtheorem{prop}{Proposition}
\newtheorem{lemma}{Lemma}
\newtheorem{cor}{Corollary}

\theoremstyle{definition}
\newtheorem{rem}{Remark}
\newtheorem{exmp}{Example}

\begin{document}

\title[Quadrature and Quantization]%
{Infinite-Dimensional Quadrature and Quantization}

\author[]
{Steffen Dereich}
\address{Steffen Dereich\\
Institut f\"ur Mathematik\\
MA 7-5, Fakult\"at II\\
Technische Universit\"at Berlin\\
Stra\ss e des 17.\ Juni 136\\
10623 Berlin\\
Germany}
\email{dereich@math.tu-berlin.de}

\author[]
{Thomas M\"uller-Gronbach}
\address{Thomas M\"uller-Gronbach\\
Institut f\"ur Mathematische Stochastik\\
Fakult\"at f\"ur Mathematik\\
Universit\"at Magdeburg\\
Postfach 4120\\
39016 Magdeburg\\
Germany}
\email{gronbach@mail.math.uni-magdeburg.de}

\author[]
{Klaus Ritter}
\address{Klaus Ritter\\
Fachbereich Mathematik\\
Technische Universit\"at Darmstadt\\
Schlo\ss gartenstra\ss e 7\\
64289 Darmstadt\\
Germany}
\email{ritter@mathematik.tu-darmstadt.de}

\date{January 6, 2006}

\begin{abstract}
We study numerical integration of Lipschitz functionals on a 
Banach space by means of deterministic and randomized (Monte Carlo)
algorithms. This quadrature problem is shown to be closely
related to the problem of quantization of the underlying 
probability measure.
In addition to the general setting we analyze in particular
integration w.r.t.\ Gaussian measures and distributions
of diffusion processes.
We derive lower bounds for the worst case
error of every algorithm in terms of its computational
cost, and we present matching upper bounds, up to
logarithms, and corresponding almost optimal algorithms.
As auxiliary results 
we determine the asymptotic behaviour of quantization numbers
and Kolmogorov widths for diffusion processes. 
\end{abstract}

\maketitle

\section{Introduction}\label{s1}
\setcounter{equation}{0}

Let $\mu$ be a Borel probability measure on a Banach space
$(\XX,\|\cdot\|)$ such that
\[
\int_{\XX} \|x\| \, \mu (dx) < \infty.
\]
Moreover, let $F$ denote the class of all Lipschitz continuous 
functionals $f:\XX\to\R$ with Lipschitz constant at most one, i.e.,
\[
\phantom{\qquad \quad x,y \in \XX.}
|f(x) - f(y) | \leq \| x - y \|,
\qquad \quad x,y \in \XX.
\]
We wish to compute
\[
S(f) = \int_{\XX} f(x)\, \mu(dx)
\]
for $f \in F$ by means of deterministic or randomized (Monte Carlo)
algorithms that use the values $f(x)$
of the functional $f$ at a finite number of 
sequentially (adaptively) chosen points $x \in \XX$.
We present a worst case analysis, and we optimally relate the error
and the cost of algorithms.

The classical instance of this quadrature problem
is given by $\XX = \R^d$ and $\mu$ being
the uniform distribution on $[0,1]^d$, say, or the $d$-dimensional standard
normal distribution. See, e.g., Novak (1988) and Wasilkowski, 
Wo\'zniakowski (2001) for results and references. 
We are mainly interested in infinite-dimensional spaces $\XX$,
and in particular we study Gaussian measures $\mu$
and distributions $\mu$ of diffusion processes,
see also Wasilkowski, Wo\'zniakowski (1996) and
Pag\`es, Printems (2004).
Infinite-dimensional quadrature is applied, e.g., 
in mathematical finance and quantum physics, and
moreover it is used as a computational tool to solve
parabolic or elliptic partial differential equations.

The appropriate framework for the analysis of finite-
and infinite-dimensional quadrature problems is
provided by the real-number model of computation.
Informally, a real-number algorithm is like a C-program that
carries out exact computations with real numbers.
Furthermore, a perfect generator for random numbers from $[0,1]$
is available, and algorithms have access to the functionals
$f \in F$ via an oracle (subroutine) that provides values $f(x)$
for points $x$ from a finite-dimensional subspace $\XX_0 \subset \XX$.
The subspace may be chosen arbitrarily, but it is fixed for
a specific algorithm. If, for instance, $\mu$ is the Wiener measure
on $\XX = C([0,1])$ or, more generally, the distribution of
a diffusion process, then spaces $\XX_0$ of piecewise
linear functions are frequently used in computational practice.
The cost of an oracle call for $f(x)$ is given by the
dimension of the corresponding subspace $\XX_0$, while 
real number operations, evaluations of elementary functions,
and calls of the random number generator are performed 
at cost one. Furthermore, in the case of a diffusion process,
function values of its drift and diffusion coefficients are
provided at cost one, too. 

By $\eps_N^\deter$ and $\eps_N^\ran$  
we denote the smallest worst
case error that can be achieved by any deterministic or randomized
algorithm, resp., whose computational cost is bounded
by $N$. We wish to determine the asymptotic behaviour of
the minimal errors $\eps_N^\deter$ and $\eps_N^\ran$ and to
find algorithms with cost close to $N$ and error close to
the corresponding minimal error.
We write $a_N \preceq b_N$ for sequences
of positive real numbers $a_N$ and $b_N$ if $\sup_{N \in \N}
a_N/b_N < \infty$. Moreover, $a_N \asymp b_N$ means
$a_N \preceq b_N$ and $b_N \preceq a_N$. 

Our main results are the following upper and lower bounds
for the minimal errors in the diffusion case and in the 
Gaussian case.

Suppose that $\mu$ is the distribution of an $m$-dimensional diffusion
process on the space $\XX = C([0,1],\R^m)$, equipped with the supremum
norm. Under mild assumptions on the 
drift and diffusion coefficients 
the minimal errors $\eps_N^\deter$ and $\eps_N^\ran$  satisfy
\[
\eps_N^\deter \succeq (\ln N)^{-1/2}
\]
see Theorem \ref{theox-1} and Proposition \ref{pqw},
and
\[
N^{-1/4} \cdot (\ln N)^{-3/4} \preceq \eps_N^\ran \preceq
N^{-1/4} \cdot (\ln N)^{ 1/4},
\]
see Theorems \ref{t7} and \ref{t8}.

We conclude that the quadrature problem for diffusion processes
is intractable
by means of deterministic algorithms, since $(\ln N)^{-1/2}$
tends to zero too slowly,
but randomization helps substantially.
The upper bound for $\eps_N^\ran$ is achieved by a suitably adjusted
weak Euler scheme. In view of the lower bound,
this algorithm is optimal, up to a multiple of at most $\ln N$,
in the class of all randomized algorithms.

Suppose that $\mu$ is a zero mean Gaussian measure on a 
separable Banach space $\XX$, whose small ball function 
\[
\varphi (\eps) = -\ln \mu(\{ x \in \XX : \|x\|  \leq \eps\})
\]
satisfies 
\[
\varphi (\eps) \asymp \eps^{-\alpha} \cdot (\ln \eps^{-1})^\beta
\]
for some constants $\alpha > 0$ and $\beta \in \R$
as $\eps$ tends to zero. 
This asymptotic behaviour 
typically holds for
Gaussian measures on infinite-dimensional spaces, see, e.g., the 
review article by Li, Shao (2001).
Consider, for instance, the distribution $\mu$ of the fractional 
Brownian motion with Hurst parameter $H \in \left]0,1\right[$ 
on the space
$\XX=C([0,1])$ or $\XX = L_p([0,1])$ with $p \in \left[1,\infty\right[$.
Then $\alpha = 1/H$ and $\beta = 0$.
A non-zero constant $\beta$ appears, for example,
in case of $\mu$ being the distribution
of the $d$-dimensional
Brownian sheet on the space $\XX = L_2([0,1]^d)$. Then $\alpha=2$
and $\beta = 2(d-1)$, see Cs\'aki (1984) and Fill, Torcaso (2004).

Given the above small ball asymptotics,
the minimal error $\eps_N^\deter$ satisfies
\[
\eps_N^\deter \succeq 
(\ln N)^{-1/\alpha} \cdot (\ln \ln N)^{\beta/\alpha},
\]
see Theorem \ref{theox-1} and Proposition \ref{pqwg}. Furthermore, for
the minimal error $\eps_N^\ran$ we have
\[
\eps_N^\ran \preceq 
N^{- 1/(2+\alpha)} \cdot (\ln N)^{(\alpha+\beta)/(2+\alpha)},
\]
see Theorem \ref{galg}, as well as
\[
\limsup_{N \to \infty}\,
\eps_N^\ran \cdot N^{1/(2+\alpha)} \cdot 
(\ln N)^{(2+2\alpha-\alpha\beta)/(\alpha(2+\alpha))} \cdot
(\ln \ln N)^{- 2\beta/(\alpha(2+\alpha))} 
>0,
\]
see Theorem \ref{tlg}. Note that the bounds for $\eps_N^\ran$ 
only differ by powers of $\ln N$ and $\ln \ln N$
for an infinite sequence of integers $N$.

We observe intractability of the quadrature problem 
for Gaussian measures by means of deterministic methods, 
and randomization helps substantially.
The upper bound for $\eps_N^\ran$ is achieved by the classical Monte
Carlo algorithm based on a normal distribution on a properly chosen
subspace $\XX_0 \subset \XX$.

We briefly outline the content of the paper.
For the analysis of the quadrature problem
we establish general relations
of the minimal errors $\eps_N^\deter$ and $\eps_N^\ran$
to quantization numbers
and average Kolmogorov and linear widths of probability measures
on Banach spaces. 
See, e.g., 
Creutzig (2002), Dereich (2003),
Dereich \emph{et al.} (2003),
Graf, Luschgy (2000),
Luschgy, Pag\`es (2004),
and Ritter (2000)
for results and references concerning the latter quantities.

In Section \ref{s2} we only use the fact 
that algorithms evaluate the functionals $f \in F$
at a finite number of points $x \in \XX$.
The minimal errors $e^\deter_n$ and $e^\ran_n$ are defined
as the smallest worst case error that can be achieved by any 
algorithm that uses $n$ functional evaluations (on average).
These minimal errors turn out to be closely related to
the quantization numbers $\qr_n$, which
are defined as a distance of the measure $\mu$ to the class of all
discrete probability measure on $\XX$ with support of size $n$.
More precisely, we have
\[
e_n^\deter = \qe_n,
\]
see Theorem \ref{theox-1}, and
\[
n^{1/2}\cdot \sup_{m \geq 4n} (\qe_{m-1}-\qe_m) \preceq
e_n^{\ran} \preceq n^{-1/2}\cdot q^{(2)}_{\lfloor n/2 \rfloor}, 
\]
see Theorems \ref{t2} and \ref{theo0727-1}. 
The latter estimate yields the well-known result $e_n^{\ran} \asymp
n^{-1/2-1/d}$ in the finite-dimensional case $\XX = \R^d$.

In Section \ref{s4} we examine the computational cost
of algorithms more closely, and 
we take into account that the functionals
$f \in F$ may only be evaluated at points $x$ from 
finite-dimensional subspaces $\XX_0 \subset \XX$. 
The latter restriction leads to the consideration
of average Kolmogorov widths $d_k^{(p)}$ of the measure $\mu$,
which are defined as average errors of best approximation
by means of optimally chosen $k$-dimensional subspaces $\XX_0 
\subset \XX$.
It turns out that
\[
\eps_N^\ran \geq 
\inf_{k \cdot n \leq N} \max (e_n^\ran,\de_k) 
\]
for every measure $\mu$, see Proposition \ref{p88}, which is the key 
tool to derive the lower bounds for randomized algorithms.. 

In Sections \ref{s5} and \ref{s6}
we study diffusion processes and Gaussian measures,
resp., and we apply the general results from Sections
\ref{s2} and \ref{s4}. As auxiliary results 
we determine the asymptotic behaviour of the quantization numbers
and the Kolmogorov widths in the diffusion case,
see Proposition \ref{pqw} and Remark \ref{r20}.

\section{Quadrature of Lipschitz Functionals and Quantization}\label{s2}

At first we disregard the details of the real number model.
We only take into account that algorithms may only evaluate 
the functionals $f \in F$ at a finite number of sequentially
chosen points in the Banach space $\XX$
in order to approximate the integrals $S(f)$.

\subsection{Basic Definitions}

Any deterministic sequential evaluation is formally defined by a point
\[
x_1\in\XX
\]
and a sequence of mappings 
\[
\phantom{,\qquad \ell \geq 2.}
\psi_\ell: \R^{\ell-1}\to \XX,\qquad \ell \geq 2.
\]
For every $f\in F$ the evaluation starts at the point $x_1$,
and the mappings $\psi_\ell$ determine the subsequent evaluation 
points. More precisely, after $n$ steps the functional values
\[
y_1 = f(x_1)
\]
and 
\[
\phantom{\qquad \ell= 2,\ldots,n,}
y_\ell = f(\psi_\ell(y_1,\ldots,y_{\ell-1})), 
\qquad \ell =2,\ldots,n,
\]
are known. 
A decision to stop or to further evaluate $f$ is made after
each step. This is formally described by a sequence of mappings
\[
\phantom{\qquad \ell \geq 1,}
\tau_\ell: \R^{\ell}\to \{0,1\},
\qquad \ell \geq 1,
\]
and the total number $n(f)$ of evaluations is given by
\[
n(f) = \min \{ \ell \geq 1 : \tau_{\ell}(y_1,\dots,y_\ell) = 1\},
\]
which is finite for every $f \in F$ by assumption.
Finally, an approximation 
\[
\Sh(f) = \phi_{n(f)}(y_1,\ldots,y_{n(f)})
\]
to $S(f)$ is defined by a sequence of mappings
\[
\phantom{\qquad \ell \geq 1.}
\phi_\ell :\R^\ell \to \R,
\qquad \ell \geq 1.
\] 

Any such mapping $\Sh : F \to \R$ could be considered as a
deterministic algorithm, with algorithm being understood in 
a broad sense,
and the corresponding class of mappings
is denoted by $\SSS^\deter$. For convenience, we identify $\Sh$
with the point $x_1$ and the sequences of mappings 
$\psi_\ell$, $\tau_\ell$, and $\phi_\ell$.
Moreover, we write $\card(\Sh,f)$ instead of $n(f)$,
and this quantity is called the cardinality of $\Sh$ applied
to $f$. 
Note that $\SSS^\deter$ contains in  particular
all quadrature formulas
\[
\Sh (f) = \sum_{i=1}^n a_i \cdot f(x_i)
\]
with $a_i \in \R$ and $x_i \in \XX$. Here
all mappings $\psi_\ell$ and $\tau_\ell$ are constant with
$\psi_2 = x_2, \dots, \psi_n = x_n$
and 
$\tau_1 = \dots = \tau_{n-1} = 0$ while $\tau_n = 1$, 
i.e., all functionals $f \in F$ are
evaluated non-sequentially at the same set of $n$ points,
and $\phi_n$ is linear.

A randomized (or Monte Carlo) broad sense algorithm based on 
sequential evaluation 
is formally defined by a probability space
$(\Omega,\A,P)$ and a mapping 
\[
\Sh : \Omega \times F \to \R
\]
such that
\begin{itemize}
\item[(i)]
$\Sh(\omega,\cdot) \in \SSS^\deter$ for every $\omega \in \Omega$,
\item[(ii)]
$\Sh(\cdot,f)$ is measurable for every $f \in F$,
\item[(iii)]
$\omega \mapsto \card (\Sh(\omega,\cdot),f)$ 
is measurable for every $f \in F$.
\end{itemize}
We refer to Nemirovsky, Yudin (1983)
and Wasilkowski (1989) for this and an equivalent
definition of randomized algorithms.
In the sequel the random variables from (ii) and (iii) are denoted 
by $\Sh(f)$ and $\card(\Sh,f)$, respectively.

By $\SSS^\ran$ we
denote the class of all mappings $\Sh$ with properties
(i)--(iii) on any probability space.
Clearly, $\SSS^\deter \subsetneq \SSS^\ran$.
Note that $\SSS^\ran$ contains in particular the 
classical (abstract) Monte Carlo method
\begin{equation}\label{g19}
\Sh(f) = 1/n \cdot \sum_{i=1}^n f(X_i)
\end{equation}
with $X_1,\dots,X_n$ being independent and distributed
according to $\mu$.

The worst case error of $\Sh \in \SSS^\ran$ is defined by
\[
e (\Sh) = \sup_{f \in F} \left( \E | S(f) - \Sh (f) |^2 \right)^{1/2},
\]
which in particular for $\Sh \in \SSS^\deter$ reads 
\[
e (\Sh) = \sup_{f \in F} | S(f) - \Sh (f) |.
\]
The worst case cardinality of $\Sh \in \SSS^\ran$ is defined by
\[
\card (\Sh) = \sup_{f \in F} \E (  \card (\Sh,f)),
\]
which in particular for $\Sh \in \SSS^\deter$ reads 
\[
\card (\Sh) = \sup_{f \in F} \card (\Sh,f). 
\]
For simplicity we assume that $\card(\Sh) \in \N$ for
randomized algorithms, too.

Minimization of the worst case error 
among those broad sense algorithms that use at most
$n$ evaluations (on average)
leads to the definition of the $n$-th minimal errors
\[
e^\ran_{n} = \inf \{ e(\Sh) : \Sh \in \SSS^\ran,\ \card (\Sh) \leq n \}
\]
and
\[
e^\deter_{n} = \inf \{ e(\Sh) : \Sh \in \SSS^\deter,\ \card (\Sh)\leq
n \}.
\]
We add that minimal errors are key quantities in 
information-based complexity, see, e.g.,
Traub, Wasilkowski, Wo\'zniakowski (1988),
Novak (1988),
and Ritter (2000).

In Sections \ref{DA} and \ref{RA} we relate the minimal errors
to quantization numbers.
The $n$-th quantization number $\qr_n$ of order $r >0$ 
is defined as 
\[
\qr_n=\inf_{x_1,\dots,x_n\in\XX} \qr(x_1,\dots,x_n),
\]
where
$$
\qr(x_1,\dots,x_n)=
\left( \int_\XX \min_{i=1,\dots,n} \|x-x_i\|^r\,\mu(dx) \right)^{1/r},
$$
see, e.g., Graf, Luschgy (2000).
In this context a collection of points $x_1,\dots,x_n \in \XX$ is 
called a codebook for quantization of the probability measure $\mu$.
For notational convenience we let $q_n = q^{(1)}_n$
and $q = q^{(1)}$.
Note that $q_n < \infty$,
and furthermore
$\lim_{n \to \infty} q_n = 0$ if $\XX$ is separable.

\subsection{Deterministic Algorithms}\label{DA}

The quantization problem and the quadrature problem by
means of broad sense deterministic algorithms are equivalent in the
following sense.
Since $S$ is a real-valued linear
mapping on a convex and symmetric set $F$, it follows that
\begin{equation}\label{g12}
e_n^\deter = 
\inf \{ e(\Sh) : 
\text{$\Sh \in \SSS^\deter$ is a quadrature formula, 
$\card(\Sh) \leq n$}\},
\end{equation}
see Smolyak (1965), Bakhvalov (1971), and also 
Traub, Wasilkowski, Wo\'zniakowski (1988, Chap.\ 4.5).
Furthermore, for $F$ and $S$ as studied in this paper we have
\[
\inf 
\{ e(\Sh) : 
\text{$\Sh \in \SSS^\deter$ is a quadrature formula
based on $x_1,\dots,x_n$}\} = q(x_1,\dots,x_n)
\]
for every codebook $x_1,\dots,x_n \in \XX$,
see Kantorovich, Rubinstein (1958) and Gray, Neuhoff, Shields (1975).
The latter infimum is attained by the quadrature formula
\begin{equation}\label{g8}
\Sh (f) = \sum_{i=1}^n \mu(V_i) \cdot f(x_i),
\end{equation}
if $V_1,\dots,V_n$ is a corresponding
Voronoi partition of $\XX$.
An (almost) optimal codebook therefore yields an (almost) optimal quadrature 
formula \eqref{g8}, and 
the $n$-th minimal error $e_n^\deter $  
coincides with the $n$-th quantization number of order one.

\begin{theorem}\label{theox-1}
For every $n \in \N$
\[
e_n^{\deter} =  q_n. 
\]
\end{theorem}

\begin{rem}\label{r1}
There are numerous results on $e_n^\deter$ or $q_n$ for 
finite-dimensional spaces $\XX = \R^d$, see, e.g.,
Novak (1988), Graf, Luschgy (2000),
Wasilkowski, Wo\'zniakowski (1996, 2001). 

Assume $r \ge 1$. Then,
under rather
mild assumptions on $\mu$, and in particular for the 
uniform distribution on $[0,1]^d$,
the quantization numbers satisfy
\begin{equation}\label{g10}
\lim_{n \to \infty}\, q^{(r)}_n \cdot n^{1/d} = c^{(r)}
\end{equation}
with some constant $c^{(r)} = c^{(r)}(\mu,d,\|\cdot\|) >0$, 
see Graf, Luschgy (2000, Thm.\ 6.2). 
\end{rem}

\begin{rem}\label{r2}
Much less is known about $e_n^\deter$ or $q_n$ for
infinite-dimensional spaces $\XX$,
see Wasilkowski, Wo\'zniakowski (1996),
Dereich \emph{et al.} (2003),
Dereich (2003, 2004), and Luschgy, Pag\`es (2003, 2004) for results and 
references. If $\mu$ is the distribution of a diffusion process 
or a Gaussian process 
then, typically, the quantization numbers $q^{(r)}_n$
tend to zero only with logarithmic order, 
see Sections \ref{s5} and \ref{s6}.
For such processes we conclude
from Theorem \ref{theox-1} that quadrature of 
arbitrary Lipschitz functionals by means
of (broad sense) deterministic algorithms is intractable.

As an example consider the Wiener measure 
$\mu$ on $\XX = C([0,1])$ endowed
with the supremum norm. In this case 
\begin{equation}\label{g7}
\lim_{n \to \infty}\, q^{(r)}_n \cdot (\ln n)^{1/2} = c
\end{equation}
with some constant $c > 0$, see 
Dereich, Scheutzow (2005).
\end{rem}

\subsection{Randomized Algorithms}\label{RA}

We first state an upper bound for the minimal error $e_n^\ran$
in terms of the quantization number $q_n^{(2)}$, which is a
consequence of a  well-known variance reduction
technique based on quantization, see Pag\`es, Printems (2004).
Note that 
$\lim_{n \to \infty} q_n^{(2)} = 0$ if $\XX$ is separable and 
$\int_\XX \|x\|^2 \, \mu(dx) < \infty$. Under the latter assumption
the classical Monte Carlo method 
\eqref{g19} without variance reduction
only yields errors of order $n^{-1/2}$ in all non-trivial cases.

\begin{theorem}\label{t2}
For every $n \in \N$
\[
e_{2n}^\ran \leq 2 \cdot n^{-1/2} \cdot q_n^{(2)}. 
\]
\end{theorem}

\begin{proof}
Consider a codebook $x_1,\dots,x_n \in \XX$ as well as a corresponding 
Voronoi partition $V_1,\ldots,V_n$ of $\XX$.
For $f\in F$ let $J(f)$ denote the interpolation of $f$
at the points $x_i$ that is constant on the 
corresponding cells $V_i$, i.e.,
\[
J (f) = \sum_{i=1}^n f(x_i) \cdot 1_{V_i}.
\]
The deterministic broad sense algorithm \eqref{g8} 
approximates $S(f)$ by $S(J(f))$. Define a broad sense randomized algorithm 
$\Sh \in \SSS^\ran$ with $\card (\Sh) \leq 2n$ by
\begin{equation}\label{g14}
\Sh (f) =
S(J(f)) + 1/n \cdot \sum_{i=1}^n (f - J(f))(X_i)
\end{equation}
with $X_1,\dots,X_n$ being independent and distributed according
to $\mu$. 
Hence the non-deter\-ministic part of $\Sh$
consists of applying the classical Monte Carlo method \eqref{g19}
to $\widetilde{f} = f - J(f)$. It follows that
\[
e(\Sh) = n^{-1/2} \cdot
\sup_{f \in F} \left( \int_\XX \left( \widetilde{f} (x) 
- S(\widetilde{f})\right)^2 \, \mu (dx) \right)^{1/2}.
\]
Since $|\widetilde{f}(x)| \leq \min_{i=1,\dots,n} \|x - x_i\|$, we
obtain
\begin{align*}
\left( \int_\XX \left( \widetilde{f} (x) - S(\widetilde{f})\right)^2 
\, \mu (dx) \right)^{1/2} 
&\leq 
\left( \int_\XX \widetilde{f}^2 (x)  \, \mu (dx) \right)^{1/2} +
|S(\widetilde{f})|\\
& \leq q^{(2)}(x_1,\dots,x_n) + q^{(1)}(x_1,\dots,x_n) \\
& \leq 2 \cdot q^{(2)}(x_1,\dots,x_n),
\end{align*}
which completes the proof.
\end{proof}

We now turn to lower bounds for (broad sense) randomized algorithms.
In this setting a result analogous to \eqref{g12} is not 
available in general, and therefore
considerations cannot a priori be restricted 
to randomized quadrature formulas.
We use the following tool, 
which is due to Bakhvalov (1959) and
Novak (1988) and which holds for integration problems in general,
see Novak (1988, Sec. 2.2.10).

\begin{prop}\label{theox-2}
Let $m \geq 4n$, 
and suppose there are
functionals $f_1, \dots, f_m : \XX \to \R$ such that
\begin{equation}\label{g3}
\{ x \in \XX : f_i(x) \neq 0\} \cap  \{ x \in \XX : f_j(x) \neq 0\}
= \emptyset
\end{equation}
for all $i \neq j$ and
\begin{equation}\label{g4}
\sum_{i=1}^m \delta_i \cdot f_i \in F
\end{equation}
for all $ \delta_1, \dots , \delta_m \in \{ \pm 1\}$.
Then
\[
e_{n}^\ran \geq \tfrac{1}{4} \cdot
n^{1/2} \cdot \min_{i=1,\dots,m} S(f_i).
\]
\end{prop}

A proper choice of the functionals $f_i$ in Proposition \ref{theox-2}
yields a lower bound for the minimal
error $e_n^\ran$ in terms of consecutive differences
of quantization numbers.

\begin{theorem}\label{theo0727-1}
For every $n\in\N$ 
\begin{align*}
e_n^{\ran}\ge 
\tfrac{1}{8} \cdot n^{1/2}\cdot \sup_{m \geq 4n} (q_{m-1}-q_m). 
\end{align*}
\end{theorem}

\begin{proof}
For $\eps\in \left]0,1\right[$ and $m \geq 4n$
  choose $x_1,\dots,x_{m}\in\XX$ with
\begin{equation}\label{g5}
  q(x_1,\dots,x_{m})\le \eps \cdot q_{m-1} +
(1-\eps) \cdot q_{m} + \eps,
\end{equation}
and consider the functionals
$$
\nice{f_i(x)=
\tfrac{1}{2} \cdot
\max(0,\min_{j\not=i} \|x-x_j\| -\|x-x_i\|)}{\qquad i=1,\dots,m}.
$$
Clearly \eqref{g3} is satisfied and $f_1,\dots,f_m \in F$.
Consequently \eqref{g4} holds, too. 

We claim that
\begin{equation}\label{g2}
S(f_i) \ge \frac{1-\eps}2 \cdot (q_{m-1}-q_m) - \eps.
\end{equation}
It suffices to prove the statement for $i=m$. To this end
     consider a Voronoi partition $V_1,\ldots,V_m$
corresponding to $x_1,\dots,x_m$, and let 
   $U_1,\ldots,U_{m-1}$  be a Voronoi partition 
corresponding to $x_1,\dots,x_{m-1}$. 
    If $j \leq m-1$ and $x\in V_m\cap U_j$ then
$$
f_m(x)= \tfrac{1}{2} \cdot \left( \|x-x_j\|- \|x-x_m\|\right).
$$
Hence
\begin{align*}
\int_\XX f_m(x) \,\mu(dx)&= \int_{V_m} f_m(x)\,\mu(dx)\\
&= \tfrac{1}{2} \cdot
\sum_{j=1}^{m-1} \int_{V_m\cap U_j} \|x-x_j\|\,\mu(dx)-
\tfrac{1}{2} \cdot \int_{V_m} \|x-x_m\|\,\mu(dx)\\
&=\tfrac{1}{2} \cdot \sum_{j=1}^{m-1} 
\int_{(V_m\cap U_j)\cup V_j} \|x-x_j\|\,\mu(dx) 
-\tfrac{1}{2} \cdot \sum_{j=1}^m \int_{V_j} \|x-x_j\|\,\mu(dx).
\end{align*}
Note that the sets $(V_m\cap U_j)\cup V_j$ with
$j \leq m-1$ form a partition of $\XX$ as well, and   every
$x\in(V_m\cap U_j)\cup V_j$ satisfies
$$
\min_{k=1,\dots,m-1} \|x-x_k\|=\|x-x_j\|.
$$
   Thus
$$
S(f_m)=\tfrac{1}{2} \cdot \left(
   q(x_1,\dots,x_{m-1})- q(x_1,\dots,x_m)\right)
$$
  and \eqref{g2} follows from \eqref{g5}.

It remains to apply Proposition \ref{theox-2} and to let $\eps$
tend to zero.
\end{proof}

The following consequence of Theorem \ref{theo0727-1} is useful, 
in particular, for finite-dimensional spaces $\XX$.

\begin{cor}\label{cor1}
If the sequence $(q_n)_{n \in \N}$ is regularly 
varying with index $-\alpha<0$ then
\[
\liminf_{n \to \infty}\,
e_n^{\ran} \cdot n^{1/2} / q_n \geq \frac{\alpha}{2^{5+2\alpha}}.
\]
\end{cor}

\begin{proof}
Put
\[
g(m)=\sup_{\ell \ge m} \left(q_{\ell-1}-q_{\ell}\right)
\]
for $m\in\N\setminus\{1\}$ and let $\kappa>1$. Clearly, 
\begin{align*}
g(m)\ge 
\frac {q_m-q_{\lceil\kappa m\rceil}}{\lceil\kappa m\rceil- m}=
q_m \cdot 
\frac{1-q_{\lceil\kappa m\rceil}/q_m}{\lceil\kappa m\rceil- m}.
\end{align*}
Since 
$\lim_{m\to\infty} q_{\lceil\kappa m\rceil}/q_m=\kappa^{-\alpha}$ 
it follows that
$$
\liminf_{m \to \infty} \, g(m) \cdot m / q_m
\geq \frac {1-\kappa^{-\alpha}}{\kappa-1}.
$$
Letting $\kappa$ tend to one yields
$$
\liminf_{m \to \infty} \, g(m) \cdot m / q_m \geq \alpha.
$$
Combining the latter estimate and Theorem \ref{theo0727-1} completes 
the proof.
\end{proof}

\begin{rem}\label{r4}
Suppose that the quantization numbers satisfy \eqref{g10},
which typically holds in the finite-dimensional case
$\XX=\R^d$, see Remark \ref{r1}. 
Then Corollary \ref{cor1} is applicable with $\alpha = 1/d$,
and we obtain
\[
\liminf_{n \to \infty} \,
e_n^{\ran} \cdot n^{1/2+1/d}  
\geq \frac{c^{(1)}}{d \cdot 2^{5+2/d}}.
\]
A matching upper bound is provided by Theorem \ref{t2}, so that
we end up with the well-known fact 
\[
e_n^\ran \asymp n^{-1/2-1/d},
\] 
see Novak (1988, Sec.\ 2.2.6) for the case of the
uniform distribution $\mu$ on $[0,1]^d$.
\end{rem}

{}From the previous remark we conclude that, 
up to multiplicative constants, neither the upper bound    
in Theorem \ref{t2} nor the lower bound in Theorem \ref{theo0727-1} 
can be improved in general.

Corollary \ref{cor1} is not applicable,
if the quantization numbers are slowly varying, cf.\ Remark \ref{r2}.
Instead, one may use the following result.

\begin{cor}\label{c2}
Let $f:\left[0,\infty\right[ \to \left]0,\infty\right[$ be a 
convex and differentiable function. If
$$
\limsup_{n \to \infty} \, q_n / f(n) \geq 1 
$$
and 
\[
\lim_{n \to \infty} q_n = 0,
\]
then 
$$
\limsup_{n\to\infty} \, e_n^\ran /  \left(n^{1/2}  
\cdot |f^\prime|(4n+3) \right) \ge 1/8.
$$
\end{cor}

\begin{proof}
Fix $\eps\in\left]0,1\right[$. By assumption
$$
q_{m-1}\ge (1-\eps) \cdot f(m-1)
= (1-\eps) \cdot \int_{m-1}^\infty -f^\prime(s)\,ds
$$
holds for infinitely many integers $m$.
Since $q_{m-1}=\sum_{k=m}^\infty (q_{k-1}-q_{k})$, we also have
$$
q_{m-1}-q_{m} \ge (1-\eps) \cdot \int_{m-1}^{m} -f^\prime(s)\,ds 
\ge - (1-\eps) \cdot  f^\prime(m)
$$
infinitely often.
To every such $m$ we associate $n=\lfloor m/4 \rfloor$. 
Then $m\in[4n,4n+3]$ and Theorem \ref{theo0727-1} implies
$$
e_n^{\ran} \ge -(1-\eps)/ 8  \cdot n^{1/2} \cdot f^\prime(4n+3).
$$ 
Letting $\eps$ tend to zero finishes the proof.
\end{proof}

\begin{rem}\label{r5}
Suppose that the quantization numbers satisfy 
\[
q_n^{(r)} \asymp (\ln n)^{-1/2},
\]
which typically holds for diffusion processes, see
Proposition \ref{pqw}, 
and in particular for the Wiener measure, see Remark \ref{r2}.
Then Corollary \ref{c2} is applicable with 
$f(t)=c\cdot (\ln t)^{-1/2}$ for some constant $c>0$,
and we obtain
\begin{equation}\label{g15}
\limsup_{n \to \infty} \,
e_n^{\ran} \cdot n^{1/2} \cdot (\ln n)^{3/2} > 0.
\end{equation}
On the other hand,
\begin{equation}\label{g16}
\limsup_{n \to \infty} \,
e_n^{\ran} \cdot n^{1/2} \cdot (\ln n)^{1/2} < \infty
\end{equation}
by Theorem \ref{t2}. 
This upper bound is achieved by a sequence
of comparatively simple broad sense randomized algorithms, see \eqref{g14},
which are far superior to any sequence of
(broad sense) deterministic algorithms, see 
Theorem \ref{theox-1}.
Moreover, upper and lower bounds do not differ much for an infinite
sequence of integers $n$. We add that, for a large class of
diffusion processes, inequality \eqref{g15} holds true 
with 
limes superior replaced by limes inferior,
see Proposition \ref{pran}.
\end{rem}

\section{Finite-dimensional Sampling and Kolmogorov Widths}
\label{s4}

So far we have studied broad sense algorithms $\Sh \in \SSS^\ran$,
and we have expressed the quality of such an algorithm
in terms of its error $e(\Sh)$ and its 
cardinality $\card (\Sh)$.
The cardinality 
serves as a crude measure of the
cost of $\Sh$, if one assumes that
evaluation of functionals $f \in F$ is possible at any point 
$x \in \XX$ at cost one and if all further operations
are not taken into account.
Moreover, by definition of $\SSS^\ran$, a broad sense randomized
algorithm may use perfect generators for random elements 
according to any Borel probability measure on $\XX$, in particular 
according to $\mu$. These assumptions are rather unrealistic
and do not correspond to a reasonable model of computation,
and the practical relevance of algorithms like 
\eqref{g14} and upper bounds like \eqref{g16}
seems to be doubtful.
We stress that this point of view concerns
lower bounds like \eqref{g15} only in the sense that they
are `too weak'. 

It is more appropriate to take the real number model of computation
as a basis for quadrature problems.
See Traub, Wasilkowski, Wo\'zniakowski (1988)
and Novak (1995) for the definition of this model.
Informally, a real number algorithm is like a C-program that 
carries out exact computations with real numbers. Furthermore, a 
perfect generator for random numbers from $[0,1]$
as well as elementary functions like $\exp$, $\ln$, etc.
are available. 
We think that these assumptions are present at least implicitly
in most of the work dealing with quadrature problems. 
Algorithms have access to the functionals $f \in F$
via an oracle (subroutine)
that provides values $f(x)$ for points $x$ from a 
finite-dimensional subspace $\XX_0 \subset
\XX$. The subspace may be chosen
arbitrarily but it is fixed for a specific algorithm, 
and the cost for each oracle call is proportional
to the dimension of $\XX_0$.

\begin{exmp}
Consider  
the distribution $\mu$ of a diffusion process $X$ with values in 
$\XX =C([0,1],\R^m)$. Let $\Xh^{(k)}$ denote the 
Euler scheme with uniform step-size $1/(k-1)$ and piecewise 
linear interpolation, and define 
the classical Euler Monte Carlo algorithm
$\Sh^{(k)}_n$ 
by
\begin{equation}\label{mceuler}
\Sh^{(k)}_{n}(f) = 1/n\cdot\sum_{i=1}^n f(\Xh^{(k)}_i)
\end{equation} 
with independent copies 
$\Xh^{(k)}_1,\ldots,\Xh^{(k)}_n$ of $\Xh^{(k)}$.
This algorithm uses an oracle for the $k$-dimensional subspace of 
piecewise linear functions with breakpoints at $\ell/(k-1)$.
Moreover, only random numbers from $[0,1]$ are needed for the 
computation of $\Sh^{(k)}_{n}(f)$.
\end{exmp}

For simplicity we assume that the cost of an oracle call for 
functional evaluation coincides with the dimension $k$ of the 
corresponding subspace
$\XX_0$ and that real number operations as well as calls of the random
number generator and evaluations of elementary functions
are performed at cost one. Furthermore, in case
of $\mu$ being the distribution of a diffusion process, function values
of its drift and diffusion coefficients are provided at
cost one, too.
Then the total cost of a computation is given by 
$k$ times the number of oracle calls for functional evaluation plus
the total number of real number operations, calls of the random number 
generator, evaluations of elementary functions, and, 
eventually, function evaluations of drift and diffusion
coefficients.

For randomized algorithms $\Sh$ the computational cost is a random
variable, which also may depend on the integrand $f \in F$.
Analogously to $\card (\Sh)$ we therefore define 
$\cost(\Sh)$,
the worst case cost of $\Sh$, by its maximal expected cost over
the class $F$.

\begin{rem}\label{r6}
For the classical Euler Monte Carlo algorithm we have
\[
\cost (\Sh^{(k)}_{n}) \asymp k \cdot n,
\]
i.e., the cost is proportional to the product 
of the dimension of the subspace 
and the number of oracle calls for functional evaluation.
Equivalently, the cost is proportional to the product of the
number of time steps and the number of repetitions.
\end{rem}

Analogously to $e_n^\ran$ we introduce the $N$-th minimal error
\[
\eps_N^\ran = \inf \{ e(\Sh) : \text{$\Sh$ randomized algorithm with
$\cost (\Sh) \leq N$} \}
\]
in the real number model. By just counting the number of oracle
calls we get $\eps_N^\ran \geq e_N^\ran$.
To derive a lower bound for $\eps_N^\ran$
that also takes into account the
dimension of the subspaces $\XX_0$ we study the 
the $k$-th average Kolmogorov width of order $p>0$
\[
d_k^{(p)} = 
\inf\biggl\{\left(\int_\XX \dist^p(x,\XX_0)\,\mu(dx)\right)^{1/p}:
\dim(\XX_0)=k\biggr\}
\]
for the measure $\mu$.
For notational convenience we let $d_k = d_k^{(1)}$.
See, e.g., Ritter (2000, Sec.\ VII.2.5) and Creutzig (2002) 
for results and references.

The following lower bound corresponds to the extremal cases, where
either the dimension $k$ of the subspace
or the number $n$ of evaluations may be arbitrarily large.

\begin{prop}\label{p88}
For every $N \in \N$
\[
\eps_N^\ran \geq 
\inf_{k \cdot n \leq N} \max (e_n^\ran,d_k).
\] 
\end{prop}

\begin{proof}
Consider any randomized algorithm $\Sh$ with $\cost (\Sh) \leq N$,
and assume that its oracle is based on a $k$-dimensional subspace
$\XX_0 \subset \XX$. 
Define a functional $f_0 \in F$ by
$f_0=\dist(\cdot,\XX_0)$. Since $\Sh$ evaluates 
$f_0$ only at
points from $\XX_0$ we have
\[
\Sh(f_0) = \Sh(-f_0),
\]
and consequently
\[
e(\Sh) \ge \tfrac{1}{2} \cdot
 \left( \left( \E | S(f_0) - \Sh (f_0) |^2 \right)^{1/2}
 +  \left( \E | S(-f_0) - \Sh (-f_0) |^2 \right)^{1/2}\right)
 \ge S(f_0). 
\]
Hence
\[
e(\Sh) \ge \int_\XX \dist(x,\XX_0)\,\mu(dx) \geq d_k.
\]
On the other hand, put $n = \card (\Sh)$ to obtain
\[
e(\Sh) \geq e_n^\ran.
\]
We conclude that $e(\Sh) \geq \max (e_n^\ran,d_k)$ for
some $k,n \in \N$ such that $k \cdot n \leq N$.
\end{proof}

\section{Randomized Algorithms for Diffusion Processes}\label{s5}

In this section we consider the distribution $\mu$
of an $m$-dimensional
diffusion process $X$ on the space $C([0,1],\R^m)$, equipped
with the supremum norm. More precisely, $X$ is given by
\begin{equation}\label{g20}
\begin{aligned}
d X_t & = a(X_t)\,dt + b(X_t)\,dW_t,\\
X_0 & = u_0 \in \R^m
\end{aligned}
\end{equation}
for $t \in [0,1]$
with an $m$-dimensional Brownian motion $W$, and we assume that
the following conditions are satisfied:
\begin{itemize}
\item [(i)]$a:\R^m\to \R^m$ is  Lipschitz continuous  
\item [(ii)] $b:\R^m\to \R^{m\times m}$ has bounded first and second 
order partial derivatives and is of class $C^\infty$ in some 
neighborhood of $u_0$
\item[(iii)] $\det b(u_0)\not=0$
\end{itemize}

We first present bounds for the quantization numbers and the
Kolmogorov widths, see also Remark \ref{r20}.
The corresponding proofs are postponed to 
Section \ref{ss1}.

\begin{prop}\label{pqw} 
The quantization numbers $q^{(r)}_n$ satisfy
\[
q^{(r)}_n \asymp (\ln n)^{-1/2}
\]
for every $r>0$. The average Kolmogorov widths $d^{(p)}_k$ satisfy
\[
d_k^{(p)} \asymp k^{-1/2}
\]
for every $p>0$.
\end{prop}

The asymptotic behavior of the
quantization numbers stated in Proposition \ref{pqw}
is partially known. Luschgy and Pag\`es (2003) study scalar
stochastic differential equations under suitable growth and smoothness 
conditions.
In this work the upper bound is established for equations with
a strictly positive diffusion coefficient
$b : [0,1] \times \R \to \R$, and a matching lower bound
is derived if $\inf_{(t,x) \in [0,1]\times \R} b(t,x) > 0$ and
$r \geq 1$. More generally, $m$-dimensional 
diffusions with a scalar diffusion coefficient $b:[0,1]\times \R^m\to\R$
are analyzed by Dereich (2004), who determines the exact asymptotic
behavior of the quantization numbers for $r \geq 1$ 
under rather mild smoothness assumptions..
The asymptotic behavior of the Kolmogorov widths 
is determined by Maiorov (1993) for the Brownian motion.

Observing Theorem \ref{theox-1} we conclude that 
quadrature of arbitrary Lipschitz functionals is intractable
by means of deterministic algorithms. 

We next present a lower bound for the minimal error $e^\ran_n$,
which improves the estimate \eqref{g15}.
See Section \ref{ss2} for the corresponding proof.

\begin{prop}\label{pran}
The minimal errors $e_n^{\ran}$ satisfy
\[
e_n^{\ran} \succeq n^{-1/2} \cdot (\ln n)^{-3/2}.
\]
\end{prop}

Propositions \ref{p88}, \ref{pqw}, and \ref{pran}
immediately yield the following lower bound.

\begin{theorem}\label{t7}
The minimal errors $\eps_N^{\ran}$ satisfy
\[
\eps_N^\ran \succeq N^{-1/4} \cdot (\ln N)^{-3/4}.
\]
\end{theorem}

Consider the Euler Monte Carlo algorithm $\Sh_n^{(k)}$ for
equation \eqref{g20} with normally
distributed increments. More precisely, 
put $\Xh^{(k)}_{i,0} = u_0$ and define
\[
\Xh^{(k)}_{i,\ell+1} = 
\Xh^{(k)}_{i,\ell} +
1/(k-1) \cdot a\bigl(\Xh^{(k)}_{i,\ell}\bigr) +
1/\sqrt{k-1} \cdot b\bigl(\Xh^{(k)}_{i,\ell}\bigr) \cdot Z_{i,\ell}
\] 
for $i=1,\dots,n$ and $\ell = 0, \dots, k-2$.
Here $(Z_{i,\ell})_{i,\ell}$ is an independent family
of $m$-dimensional standard normally distributed random vectors.
Finally, let
$\Xh^{(k)}_i$ denote the piecewise linear interpolation
of $\Xh^{(k)}_{i,0}, \dots, \Xh^{(k)}_{i,k-1}$ at the
breakpoints $\ell/(k-1)$. Then $\Sh_n^{(k)}$
is given by \eqref{mceuler}.

\begin{theorem}\label{t8}
The Euler Monte Carlo algorithm $\Sh_N = \Sh_n^{(k)}$ 
with $n = \lfloor N^{1/2} \cdot (\ln N)^{-1/2} \rfloor$ and
$k = \lfloor N^{1/2} \cdot (\ln N)^{1/2} \rfloor$ 
satisfies
\[
e(\Sh_N) \preceq N^{-1/4} \cdot (\ln N)^{1/4}
\]
and
\[
\cost (\Sh_N) \preceq N.
\]
\end{theorem}

\begin{proof}
Consider the strong Euler scheme $\Xh^{(k)}$ with step-size $1/(k-1)$
and piecewise linear interpolation for equation \eqref{g20}. 
Then
\begin{equation}\label{g22}
\E \| X - \Xh^{(k)} \|_\infty
 \leq c_1 \cdot k^{-1/2} \cdot (\ln k)^{1/2}
\end{equation}
with some constant $c_1 > 0$
that does not depend on $k$, see Faure (1992). Let $f \in F$.
Since $S(f) = \E( f(X))$ and
$\E (\Sh_N (f)) = \E( f(\Xh^{(k)}))$, we get 
\[
| S(f) - \E (\Sh_N(f)) |
\leq c_1 \cdot\left(  \ln k / k \right)^{1/2}
\]
for the bias of $\Sh_N(f)$
by means of \eqref{g22}.
Put $g = f - f(0)$ to obtain
\[
\V (\Sh_N(f)) = \V (\Sh_N (g))
 \leq 1/n \cdot \E (g^2(\Xh^{(k)})) \leq 
1/n \cdot \E (\|\Xh^{(k)}\|_\infty^2) \leq c_2 \cdot 1/n
\]
for the variance of $\Sh_N(f)$, 
where the constant $c_2 > 0$ depends neither on $k$ nor on $f$. 
We conclude that
\[
\E (S(f) - \Sh_N(f))^2 \leq \max (c_1^2,c_2) \cdot
\left(1/n + \ln k / k \right),
\]
and with the particular choice of $n$ and $k$ the asymptotic
upper bound for the error of $\Sh_N = \Sh_n^{(k)}$ follows.
The cost of $\Sh_N$ is determined in Remark \ref{r6}.
\end{proof}

Combine Theorems \ref{t7} and \ref{t8} to conclude that
the Euler Monte Carlo algorithm $\Sh_N$ is almost
optimal.

\begin{cor}\label{c3}
\[
e(\Sh_N) \preceq \eps^\ran_N \cdot \ln N.
\]
\end{cor}

\subsection{Preliminaries}

A basic idea in the proofs of Propositions \ref{pqw} and
\ref{pran} is to reduce the 
the case of an $m$-dimensional diffusion process with properties
(i)--(iii)
to the particular case of a one-dimensional Brownian
motion by means of Lipschitz transformations and stopping.

Let $X$ denote any random element
with values in some Banach space $\XX$ and consider its distribution
$\mu$ on this space. We use the notation 
\[
e^\ran_n(X,\XX) = e^\ran_n
\]
for the $n$-th minimal error of randomized algorithms,
\[
d_k^{(p)} (X,\XX) = d_k^{(p)} 
\]
for the $k$-th average Kolmogorov width of order $p$, and
\[
\qr_n (X,\XX) = \qr_n
\]
for the $n$-th quantization number of order $r$.

Consider a measurable mapping $T : \XX \to \YY$, where
$\YY$ is a Banach space, too.
The following observation is straightforward to verify.
We add that an analogous result for Kolmogorov widths is not available.

\begin{lemma}\label{l1} 
Suppose that $T$ is Lipschitz continuous 
with a Lipschitz constant $L > 0$. Then
$$
e_n^\ran(TX,\YY)\le L \cdot e_n^\ran(X,\XX)
$$
and
\[
\qr_n(TX,\YY)\le L \cdot \qr_n(X,\XX).
\]
\end{lemma}

We formulate a simplified version of a general relation
between quantization numbers and average
Kolmogorov widths, which is due to Creutzig (2002, Thm.\ 4.6.1).

\begin{lemma}\label{carl} 
For $0 < r < p$
$$
\sup_{n\le 2^\ell} \, \ln n \cdot \qr_n (X,\XX)
\preceq
\sup_{k\le \ell} \, k \cdot d_k^{(p)}(X,\XX).
$$
\end{lemma}

The following contraction principle holds for best approximation
of sums of independent and symmetric random elements.

\begin{lemma}\label{ljakob}
Let $X_1,\dots,X_k$ denote a sequence of independent and
symmetric random elements with values in $\XX$ and let $p \geq 1$. Then
\[
\E \Bigl( \dist^p\Bigl(\sum_{\ell=1}^k \lambda_\ell X_\ell, \XX_0\Bigr) 
\Bigr)
\leq
\max_{\ell=1,\dots,k} |\lambda_\ell|^p \cdot
\E \Bigl( \dist^p\Bigl(\sum_{\ell=1}^k X_\ell, \XX_0\Bigr) \Bigr)
\]
for all $\lambda_1,\dots,\lambda_k \in \R$ and
every closed linear subspace $\XX_0 \subset \XX$.
\end{lemma}

\begin{proof}
Take Rademacher variables $\eps_1,\dots,\eps_k$ such
that $\eps_1,\dots,\eps_k,X_1,\dots,X_k$ are independent, and
consider the quotient mapping $Q : \XX \to \XX / \XX_0$.
Since $(X_1,\dots,X_k)$ and $(\eps_1 X_1,\dots,\eps_k X_k)$
coincide in distribution, the same holds true for
$(Q X_1,\dots,QX_k)$ and $(\eps_1  QX_1,\dots,\eps_k Q X_k)$.
Hence
\[
\E \Bigl( \dist^p\Bigl(\sum_{\ell=1}^k \lambda_\ell X_\ell, \XX_0\Bigr)
 \Bigr)
=
\E \Bigl\| \sum_{\ell=1}^k \lambda_\ell \cdot Q X_\ell 
\Bigr\|_{\XX / \XX_0}^p
=
\E \Bigl\| 
\sum_{\ell=1}^k \lambda_\ell \eps_\ell \cdot Q X_\ell 
\Bigr\|_{\XX / \XX_0}^p.
\]
For any choice of elements $y_\ell \in \XX / \XX_0$
\[
\E \Bigl\| \sum_{\ell=1}^k \lambda_\ell \eps_\ell \cdot y_\ell 
\Bigr\|^p_{\XX/\XX_0}
\leq
\max_{\ell=1,\dots,k} |\lambda_\ell|^p \cdot
\E \Bigl\| 
\sum_{\ell=1}^k \eps_\ell \cdot y_\ell \Bigr\|_{\XX / \XX_0}^p
\]
due to Kahane's contraction principle, see
Kahane (1993, p.\ 21). Thus
\[
\E 
\Bigl\| \sum_{\ell=1}^k \lambda_\ell \eps_\ell \cdot Q X_\ell 
\Bigr\|_{\XX / \XX_0}^p
\leq
\max_{\ell=1,\dots,k} |\lambda_\ell|^p \cdot
\E
\Bigl\| 
\sum_{\ell=1}^k \eps_\ell \cdot Q X_\ell \Bigr\|_{\XX / \XX_0}^p,
\]
which completes the proof.
\end{proof}

Now we turn to the diffusion process $X$ given by \eqref{g20}.

\begin{lemma}\label{le1124-1}
There exists a neighborhood $U$ of $u_0$ and a function 
$f \in C^\infty(U)$ such that
\[
(\grad f)^* bb^* \grad f = 1.
\]
\end{lemma}

\begin{proof}
Choose a radius $r > 0$ such that $\det bb^* (u) \neq 0$
if $|u-u_0| < r$. Furthermore, take $g \in C^\infty(\R^m,\R^{m
\times m})$ with symmetric and positive definite values
such that 
\[
g(u) = (bb^*)^{-1} (u)
\]
if $|u-u_0| < r/2$ and $g(u)$ is the identity matrix if $|u-u_0| > r$.
Then $M=\R^m$ endowed with the metric tensor 
$\sum_{i,j=1}^m g_{ij}(u) \cdot du^i\otimes du^j$ is a complete
$C^\infty$-Riemannian manifold. 
Here $u^1,\dots,u^m$ are the local coordinates obtained when taking the 
identity as chart. Moreover, let $d_M$ denote the corresponding
Riemannian distance.

Choose $v_0 \in M$ such that 
$0 < |v_0 - u_0| < r/2$ and
$0 < d_M(v_0, u_0) < i_{v_0} (M)$,
where $i_{v_0}(M)$ denotes the injectivity radius at $v_0$,
see Sakai (1996, Prop.\ III.4.13).
Define
\[
U = \{ u \in M : 0 < |v_0 - u| < r/2,\ 0<d_m(v_0,u) < i_{v_0} (M)\}
\]
as well as
\[
f(u) = d_M(v_0,u)
\]
for $u \in U$.
Then $f \in C^\infty(U)$ and
$(\grad f)^* bb^* \grad f = 1$, see
Sakai (1996, Prop..\ III.4.8).
\end{proof}

In addition to 
$C = C([0,1],\R^m)$ we also consider
the Banach space  $L_1 = L_1([0,1],\R^m)$.

\begin{lemma}\label{l2}
Either let $\XX=C$ and $\YY=C([0,1],\R)$ or let 
$\XX=L_1$ and $\YY=L_1([0,1],\R)$. There exists a Lipschitz continuous 
mapping $T:\XX\to\YY$ and a stopping time 
$\tau$ with $\P (\tau>0)=1$ such that the stopped process
\[
\phantom{,\qquad t\in[0,1],}
(T X)^\tau_t = (T X)_{t \wedge \tau},\qquad t\in[0,1],
\]
is a Brownian motion stopped at time $\tau$.
\end{lemma}

\begin{proof}
Due to Lemma \ref{le1124-1} there exists a function 
$h\in C^\infty(\R^m)$ with bounded derivatives that satisfies 
\begin{equation}\label{m1}
(\grad h)^* bb^* \grad h = 1
\end{equation}
on a closed ball with radius $r>0$ around $u_0$. 
Define the stopping time
\[
\tau = \inf\{t\in[0,1]:\, |X_t-u_0|=r\}.
\]
Clearly, $\P(\tau>0)=1$.
 
In both cases cases, $\XX=C$ 
and $\XX=L_1$ we define a Lipschitz continuous mapping $T:\XX\to\YY$ by
\[
(T x)(t)= h(x(t))- h(u_0)- 
\int_0^t \Bigl((\grad h)^*a +
\tfrac{1}{2}\cdot \sum_{i,j=1}^m (bb^*)_{i,j} 
\tfrac {\partial^2}{\partial u^i\partial  u^j} h\Bigr)(x(s))\,ds.
\]
It\^{o}'s formula implies
\[
(TX)_t = \int_0^t \bigl((\grad h)^*b\bigr)(X_t)\, dW_t. 
\]
Observing \eqref{m1} we conclude that the stopped process $(TX)^\tau$ 
is a continuous martingale with quadratic variation
\[
\langle (TX)^\tau \rangle_t = \int_0^{t\wedge \tau}
\bigl((\grad h)^* bb^* \grad h\bigr) (X_s)\, ds = t\wedge \tau,
\]
which completes the proof.
\end{proof}

\begin{rem}
The assumption that the diffusion
coefficient $b$ is of class $C^\infty$ in a neighborhood of the
initial value $u_0$ can be relaxed.
For instance, in the one-dimensional case it suffices to 
assume $b \in C^1([0,1])$ with Lipschitz continuous first derivative. 
Then
$$
f(u)=\int_{u_0}^u |1/b(v)|\,dv 
$$
is well defined in a neighborhood of $u_0$, and the statement of 
Lemma \ref{l2} follows with the same proof.
\end{rem}

\subsection{Proof of Proposition \ref{pqw}}\label{ss1}

We use the contraction principle from Lemma \ref{ljakob} to establish
the upper bound for the Kolmogorov widths.

\begin{lemma}\label{h1} 
For every $p > 0$
\[
d_k^{(p)}(X,C) \preceq k^{-1/2}.
\]
\end{lemma}

\begin{proof}
Assume that $p \geq 1$ without loss of generality.
Fix $k\in\N$, put $t_\ell=\ell/k$ for $\ell=0,\ldots,k$, and consider
the corresponding Euler process $\Xb^{(k)}$ defined by $\Xb^{(k)}_0=u_0$
and 
\[
\Xb^{(k)}_t = \Xb^{(k)}_{t_\ell} + a(\Xb^{(k)}_{t_\ell})\cdot (t-t_\ell)
+ b(\Xb^{(k)}_{t_\ell})\cdot (W_t-W_{t_\ell})
\] 
for $t\in[t_\ell,t_{\ell+1}]$. We have
\[
\E\|X-\Xb^{(k)}\|_\infty^p \preceq k^{-p/2},
\]
see Bouleau, L\'epingle (1994, p.\ 276), and therefore
\[
d_k^{(p)}(X,C) \preceq k^{-1/2}+ d_k^{(p)}(\Xb^{(k)},C).
\]

Let $\Wt^{(k)}$ denote the piecewise linear interpolation of the 
Brownian motion $W$ at the breakpoints $t_\ell$ and
define the continuous process $V^{(k)}$ by
\[
V^{(k)}_t = b(\Xb^{(k)}_{t_\ell})\cdot (W_t-\Wt^{(k)}_t)
\]
for $t\in[t_\ell,t_{\ell+1}]$. Note that $\Xb^{(k)}-V^{(k)}$
takes values in the 
$(k+1)$-dimensional subspace of piecewise linear functions
with breakpoints $t_\ell$. 
Hence
\[
d_{2k+1}^{(p)}(\Xb^{(k)},C)\le d_{k}^{(p)}(V^{(k)},C).
\]
Let $\A$ denote the $\sigma$-algebra generated by $W(t_1),\dots,W(t_k)$.
The random variables $b(\Xb^{(k)}_{t_\ell})$ are measurable with respect
to $\A$, and conditioned on $\A$ the process $W-\Wt^{(k)}$ consists of 
independent Brownian bridges on the subintervals $[t_\ell,t_{\ell+1}]$.
We apply Lemma \ref{ljakob} with $X_\ell = 1_{[t_{\ell-1},t_\ell]}
\cdot (W -\Wt^{(k)})$ to obtain
\[
d^{(p)}_{2k} (V^{(k)},C) 
\leq \bigl(\E \|b(\Xb^{(k)})\|_\infty^p \bigr)^{1/p}
\cdot d^{(p)}_{2k} (W - \Wt^{(k)},C)
\preceq d^{(p)}_k (W,C).
\]
{}From Maiorov (1993) we get $d^{(p)}_k (W,C) \asymp k^{-1/2}$.
\end{proof}

The lower bound for the quantization numbers even holds for the
space $\XX = L_1$.

\begin{lemma}\label{h2} 
For every $r > 0$
\[
q_n^{(r)}(X,L_1) \succeq (\ln n)^{-1/2}.
\]
\end{lemma}

\begin{proof}
Observe that,
due to Lemma \ref{l1} and Lemma \ref{l2},
it suffices to show that
\begin{equation}\label{eq2}
\qr_n (Y, L_1) \succeq (\ln n)^{-1/2}
\end{equation}
for every one-dimensional process $Y$ such that
\[
\phantom{,\qquad t \in [0,1],}
Y_{t \wedge \tau} = W_{t \wedge \tau},\qquad t \in [0,1], 
\]
with a stopping time $\tau$ that satisfies
$P(\tau = 0)= 0$.

To this end fix $\eps \in \left]0,1\right]$ with 
$\P(\tau \geq \eps) > 0$
and define a bounded linear operator $T : L_1 \to L_1$ by
\[
(T x)(t) = \eps^{-1/2} \cdot x(\eps \cdot t).
\]
Clearly $TW$ is a Brownian motion, too.
The quantization problem for Gaussian processes in the space $L^1$ 
is analyzed in Dereich, Scheutzow (2005).
In particular there exists a constant $\kappa>0$ such that 
\begin{equation}\label{DS1}
\lim_{n\to\infty}\, (\ln n)^{1/2}\cdot q^{(r)}_n(TW,L_1) =  \kappa
\end{equation}
for every $r>0$, see Dereich, Scheutzow (2005, Thm. 6.1).

For $n \in \N$ let $M_n \subset L_1$ denote any set of cardinality $n$,
fix $\delta \in \left]0,1\right[$, and put
\[
A_n = \{ \dist(TW,M_n) \geq (1-\delta) \cdot \qr_n (TW,L_1)\}.
\]
Due to \eqref{DS1} we can complement the sets $M_n$ to sets 
$\widetilde M_n$ of cardinality $2n$ such that
\[
\lim_{n\to\infty} \, (\ln n)^{1/2} \cdot 
\bigl( \E (\dist^{2r}(TW,\widetilde M_n))\bigr)^{1/2r}=\kappa. 
\]
as well as
\[
\lim_{n\to\infty} \, (\ln n)^{1/2} \cdot 
\bigl( \E (\dist^{r}(TW,\widetilde M_n))\bigr)^{1/r}=\kappa.
\]
Employing Lemma A.1 in Dereich, Scheutzow (2005) we conclude that
\[
\lim_{n\to\infty} \P(A_n)=1.
\]
Consequently
\begin{align*}
\E (\dist^r (TY, M_n)) &\geq 
\E (1_{\{\tau \geq \eps\}} \cdot \dist^r (TW, M_n)) \\
&\ge
(1-\delta)^r \cdot \P (\{\tau \geq \eps \} \cap A_n) \cdot
\left(\qr_n(TW,L_1)\right)^r \\
&\succeq (\ln n)^{-r/2},
\end{align*}
which yields
\[
\qr_n (TY, L_1) \succeq (\ln n)^{-1/2}.
\]
The latter bound implies \eqref{eq2} by Lemma \ref{l1}.
\end{proof}

\begin{proof}[Proof of Proposition \ref{pqw}]
In view of Lemma \ref{h1} and Lemma \ref{h2} it suffices to show
that
\begin{equation}\label{h3}
q_n^{(r)}(X,C) \preceq (\ln n)^{-1/2}
\end{equation}
and
\begin{equation}\label{h4}
d_k^{(p)}(X,L_1) \succeq k^{-1/2}.
\end{equation}

By Lemma \ref{carl} and Lemma \ref{h1} we have
\[
\ln n \cdot \qr_n (X, C) \preceq
\sup_{k\le 2\ln n} \, k \cdot d_k^{(2r)}(X, C) 
\preceq (\ln n)^{1/2},
\] 
which yields \eqref{h3}. From Lemma \ref{h1} we also get 
\begin{equation}\label{h5}
d_k^{(p)} (X,L_1) \leq c \cdot k^{-1/2}
\end{equation}
with some constant $c>0$. Moreover, by Lemma \ref{h2}, 
\[
\sup_{n\le 2^\ell} \, \ln n \cdot q_n^{(p/2)} (X,L_1)
\succeq \sup_{n\le 2^\ell} (\ln n)^{1/2} \succeq \ell^{1/2}.
\]
Consequently, by Lemma \ref{carl}
\[
\sup_{k\le \ell} \, k \cdot d_k^{(p)}(X,L_1) \geq 
\widetilde{c} \cdot \ell^{1/2}
\]
with some constant $\widetilde{c} \in \left]0,c\right[$.
Put $c = (\widetilde{c}/c)^2$.
Since
\[
\sup_{k< c\cdot \ell} \, k \cdot d_k^{(p)}(X,L_1) < 
\widetilde{c} \cdot \ell^{1/2}
\]
by \eqref{h5}, we conclude that
\[
\ell \cdot d_{\lfloor c \cdot \ell \rfloor}^{(p)}(X,L_1)
\geq \sup_{c \cdot \ell \leq k\le \ell} \, k \cdot d_k^{(p)}(X,L_1) 
\geq \widetilde{c} \cdot \ell^{1/2},
\]
which yields 
\eqref{h4}.
\end{proof}

\begin{rem}\label{r20}
According to Lemma \ref{h2} and \eqref{h4}, 
Proposition \ref{pqw} is valid, too, for $\XX=L_1$
instead of $\XX = C$.
\end{rem}

\subsection{Proof of Proposition \ref{pran}}\label{ss2}

Consider a one-dimensional Brownian motion $W$. 
Given $\ell \in\N$ and $\eps \in \left]0,1\right]$ let 
$s_i=i/ \ell \cdot \eps$ and put
\[
B^{\ell,\eps}_{i,0} = 
\{ x \in C([0,1]) : x(s_i) - x(s_{i-1}) \geq \eps^{1/2} / \ell^{3/2} \}
\]
as well as
\[
B^{\ell,\eps}_{i,1} = 
\{ x \in C([0,1]) : x(s_i) - x(s_{i-1}) < - \eps^{1/2} / \ell^{3/2} \}
\]
for $i=1,\dots,\ell$.
Moreover, define
$$
A_\alpha^{\ell,\eps} = \bigcap_{i=1}^\ell 
\{ W \in B^{\ell,\eps}_{i,\alpha_i} \}
$$
for any multi-index $\alpha\in\{0,1\}^\ell$.

\begin{lemma}\label{le0706-1}
There exists a constant $c_0 \in \left]0,1\right[$ such that 
$$
c_0 \cdot 2^{-\ell}\le \P (A_{\alpha}^{\ell,\eps})\le 2^{-\ell}
$$
for all $\ell\in\N$, $\eps \in \left]0,1\right]$, 
and $\alpha\in\{0,1\}^\ell$.
\end{lemma}

\begin{proof}
Obviously, the probability $\P(A_{\alpha}^{\ell,\eps})$ 
does not depend on $\alpha$. Hence
\begin{align*}
\P(W \in B^{\ell,\eps}_{i,0})&= 
\tfrac12 - \P (0\le W_{s_i}-W_{s_{i-1}}\le \eps^{1/2}/\ell^{3/2})\\
&=\tfrac12- \int_0^{1/\ell} (2\pi)^{-1/2} \exp(-x^2/2) \,dx\\
&\ge \tfrac12 \cdot\bigl(1-\sqrt{2/\pi} \cdot \ell^{-1}\bigr)
\end{align*}
implies 
$$
2^{\ell} \cdot \P( A_\alpha^{\ell,\eps})\ge 
\bigl(1-\sqrt{2/\pi} \cdot \ell^{-1}\bigr)^\ell.
$$
The latter bound  tends
to $\exp (-\sqrt{2/\pi})$ as $\ell$ tends to infinity, which 
completes the proof.
\end{proof}

Let $A$ be any event with $\P(A)\ge 1-c_0/2$ and put 
$$
N(\eps,\ell)= \#\{\alpha\in\{0,1\}^\ell: 
\P(A_\alpha^{\ell,\eps}\cap A)> c_0 \cdot 2^{-\ell-2}\}.
$$

\begin{lemma} \label{le0718-1} 
For all $\eps>0$ and $\ell \in \N$ 
$$
N(\eps,\ell)\ge c_1 \cdot 2^\ell,
$$
where $c_1=c_0/(4-c_0)$.
\end{lemma}

\begin{proof}
Due to Lemma \ref{le0706-1}
\begin{align*}
\P \biggl(\,\bigcup_{\alpha\in\{0,1\}^\ell} A_\alpha^{\ell,\eps}
\cap A\biggr) 
&\ge 
\P(A)+\P\biggl(\,\bigcup_{\alpha\in\{0,1\}^\ell} A_\alpha^{\ell,\eps}
\biggr)-1\\
&\ge \P(A)+ c_0-1\ge c_0/2.
\end{align*}
On the other hand, by the definition of $N(\eps,\ell)$ and 
Lemma \ref{le0706-1} 
$$
\P\biggl(\,\bigcup_{\alpha\in\{0,1\}^\ell} A_\alpha^{\ell,\eps}
\cap A\biggr)\le (2^\ell-N(\eps,\ell)) \cdot 
c_0 \cdot 2^{-\ell-2}+ N(\eps,\ell) \cdot 2^{-\ell}. 
$$
It remains to combine both estimates. 
\end{proof}

\begin{proof}[Proof of Proposition \ref{pran}]
Because of Lemma \ref{l1} and Lemma \ref{l2} it suffices to
prove that
\begin{equation}\label{eq1}
e^\ran_n(Y,C) \succeq n^{-1/2} \cdot (\ln n)^{-3/2}
\end{equation}
for every one-dimensional process $Y$
such that
\[
\phantom{,\qquad t \in [0,1],}
Y_{t \wedge \tau} = W_{t \wedge \tau},\qquad t \in [0,1], 
\]
with a stopping time $\tau$ that satisfies
$P(\tau = 0)= 0$.
To this end we use Proposition \ref{theox-2}.

Put
$$
B^\ell_\alpha = 
\bigcap_{i=1}^\ell B^{\ell,0}_{i,\alpha_i} 
$$
and define 
$f^\ell_\alpha \in F$ by
$$
f^\ell_\alpha(x) = \dist \left(x,\left(B^\ell_\alpha\right)^c\right)
$$
for $\alpha\in\{0,1\}^\ell$. Note that
\[
f^\ell_\alpha (x)\geq
\tfrac{1}{2} \cdot \min_{i=1,\dots,\ell} |x(s_i)-x(s_{i-1})|
\]
for $x\in B^\ell_\alpha$. Choose 
$\eps\in\left]0,1\right]$ with $\P(\tau\ge \eps)\ge 1-c_0/2$,
and let $A=\{\tau\ge \eps\}$. Then
\begin{align*}
S(f^\ell_\alpha)&\ge  
\tfrac{1}{2}  \cdot
\E\Bigl(1_A\cdot 1_{B^\ell_\alpha}(W) \cdot \min_{i=1,\dots,\ell}
|W_{s_i}-W_{s_{i-1}}|\Bigr)\\
&\ge \tfrac{1}{2}
\cdot  \E\Bigl(1_{A_\alpha^{\ell,\eps}\cap A} \cdot \min_{i=1,\dots,\ell}
|W_{s_i}-W_{s_{i-1}}|\Bigr)\\
&\ge \tfrac{1}{2} \cdot
\eps^{1/2}/ \ell^{3/2} \cdot \P(A_\alpha^{\ell,\eps}\cap A).
\end{align*}
Take $n = \lfloor c_1\cdot 2^{\ell-1}\rfloor$ and use 
Lemma \ref{le0718-1} to conclude that 
\[
S(f^\ell_\alpha) \succeq n^{-1} \cdot (\ln n)^{-3/2}
\]
holds uniformly for at least $2n$ multi-indices $\alpha\in\{0,1\}^\ell$.
Finally, apply Proposition \ref{theox-2} to complete the proof of \eqref{eq1}.
\end{proof}

\section{Randomized Algorithms for Gaussian Measures}\label{s6}

In this section we consider zero mean
Gaussian measures $\mu$ on separable Banach spaces $\XX$,
and throughout 
we assume that the corresponding small ball function
\[
\varphi (\eps) = -\ln \mu(\{ x \in \XX : \|x\|  \leq \eps\})
\]
satisfies
\begin{equation}\label{n1}
\varphi (\eps) \asymp \eps^{-\alpha}\cdot (\ln \eps^{-1})^{\beta}
\end{equation}
for some constants $\alpha>0$ and $\beta\in\R$ as $\eps$ tends to zero. 

\begin{rem}
Typically, \eqref{n1} holds for infinite-dimensional spaces $\XX$, 
see Li, Shao (2001).
For example, if $\mu$ is the distribution of a fractional
Brownian motion with Hurst parameter $H \in \left]0,1\right[$ on 
$\XX=C([0,1])$ or $\XX = L_p([0,1])$ for some 
$p \in \left[1,\infty\right[$, then $\alpha=1/H$
and $\beta=0$.. Moreover, 
$\alpha=1/(H-\gamma)$ and 
$\beta=0$ when $\|\cdot\|$ denotes the $\gamma$-H\"older norm. 
Similar results are known for Sobolev norms, 
see Kuelbs, Li, Shao (1995) and Li, Shao (1999).

If $\XX=C([0,1]^2)$ and $\mu$ is the distribution of the 
two-dimensional fractional Brownian sheet,
then $\alpha=1/H$ and $\beta=1+1/H$ due to Belinsky, Linde (2002). 
Moreover, for a $d$-dimensional Brownian sheet considered in 
$\XX=L_2([0,1]^d)$ one has $\alpha=2$ and $\beta=2(d-1)$, 
see Cs\'aki (1984) and Fill, Torcaso (2004).
\end{rem}

Assumption \eqref{n1} determines the asymptotic behavior of
the quantization numbers and the Kolmogorov widths,
see Dereich (2003, Thm.\ 3.1.2) and Creutzig (2002, Cor.\ 4.7.2). 
 
\begin{prop}\label{pqwg}
The quantization numbers $q^{(r)}_n$ satisfy
\[
q^{(r)}_n \asymp (\ln n)^{-1/\alpha}\cdot (\ln\ln n)^{\beta/\alpha}
\]
for every $r>0$. The average Kolmogorov widths $d^{(p)}_k$ satisfy
\[
d_k^{(p)} \asymp k^{-1/\alpha}\cdot (\ln k)^{\beta/\alpha}
\]
for every $p>0$.
\end{prop}

Hence, by Theorem \ref{theox-1}, quadrature of arbitrary Lipschitz 
functionals
by means of deterministic algorithms is intractable.
Now we turn to the analysis of randomized algorithms.

\begin{prop}\label{lower1g}
The minimal errors $e_n^\ran$ satisfy
\[
\limsup_{n \to \infty}\, 
e_n^\ran \cdot n^{1/2} \cdot (\ln n)^{1+1/\alpha}
\cdot (\ln\ln n)^{-\beta/\alpha} > 0.
\]
\end{prop}

\begin{proof}
Apply Corollary \ref{c2} with $f$ 
given by
$f(t) = c \cdot (\ln t)^{-1/\alpha}\cdot (\ln\ln t)^{\beta/\alpha}$
for $t$ sufficiently large and a suitable constant $c>0$.
\end{proof}

Proposition \ref{lower1g} provides a lower bound for the 
error of broad sense randomized algorithms in terms of the
number of functional evaluations. 
The lower bound depends on the specific properties of the Gaussian
measure only via logarithmic terms. This is no longer the case
if we relate the error of randomized algorithms to their
computational cost.

\begin{theorem}\label{tlg}
The minimal errors $\eps_N^\ran$ satisfy
\[
\limsup_{N \to \infty}\,
\eps_N^\ran \cdot N^{1/(2+\alpha)} \cdot 
(\ln N)^{(2 +2\alpha - \alpha\beta)/(\alpha(2+\alpha))} \cdot
(\ln \ln N)^{-2 \beta /(\alpha (2+\alpha))} > 0 .
\]
\end{theorem}

\begin{proof}
We combine Propositions \ref{p88}, \ref{pqwg}, and \ref{lower1g}.
Due to Proposition \ref{lower1g} there exists a constant $c>0$ and
an increasing sequence of integers $n_\ell\in\N$ such that
\[
e_{n_\ell}^\ran \ge c \cdot n_\ell^{-1/2} \cdot (\ln n_\ell)^{-1-1/\alpha}
\cdot (\ln\ln n_\ell)^{\beta/\alpha}
\]
for every $\ell\in\N$. Put
\[
N_\ell = \bigl\lfloor n_\ell^{(2+\alpha)/2} 
\cdot (\ln n_\ell)^{\alpha+\beta +1} \cdot
(\ln \ln n_\ell)^{-\beta}\bigr\rfloor,
\]
and let $n,k\in\N$ with $n\cdot k \le N_\ell$.
If $n > n_\ell$ then $k < N_\ell/n_\ell$, and Proposition \ref{pqwg} 
implies 
\begin{equation}\label{g345}
d_k \ge d_{\lfloor N_\ell/n_\ell\rfloor}
\succeq (N_\ell/n_\ell)^{-1/\alpha}
\cdot(\ln (N_\ell/n_\ell))^{\beta/\alpha}
\asymp n_\ell^{-1/2}\cdot (\ln n_\ell)^{-1-1/\alpha}
\cdot (\ln\ln n_\ell)^{\beta/\alpha}.
\end{equation}
On the other hand, if $n \le n_\ell$ then 
$e_n^\ran \ge e_{n_\ell}^\ran$. Consequently, 
by Proposition \ref{p88} and \eqref{g345}
\[
\eps_{N_\ell}^\ran \succeq n_\ell^{-1/2}\cdot (\ln n_\ell)^{-1-1/\alpha}
\cdot (\ln\ln n_\ell)^{\beta/\alpha}.
\]
Straightforward computations show 
\begin{align*}
& n_\ell^{-1/2} \cdot (\ln n_\ell)^{-1-1/\alpha}
\cdot (\ln\ln n_\ell)^{\beta/\alpha} \\
& \quad\qquad \asymp N_\ell^{-1/(2+\alpha)}\cdot 
(\ln N_\ell)^{-(2 +2\alpha - \alpha\beta)/(\alpha(2+\alpha))} \cdot
(\ln \ln N_\ell)^{2 \beta /(\alpha (2+\alpha))}, 
\end{align*}
which completes the proof.
\end{proof}

It is quite common to approximately compute the integrals
$S(f)$ with respect to Gaussian measures by sampling
from a standard normal distribution on a suitable
finite-dimensional subspace of $\XX$.
A proper choice of the subspace is suggested
by the following general result on average linear widths,
which is due to Creutzig (2002, Thm.\ 4.4.1).
There exist points $x_\ell^{(k)} \in \XX$ and
bounded linear functionals $\xi_\ell^{(k)} \in \XX^*$  
such that
\begin{equation}\label{n3}
\left( \int_{\XX} \| x - \widehat{X}^{(k)}(x) \|^2 \, 
\mu (dx) \right)^{1/2} \preceq \ln k \cdot d_k
\end{equation}
for
\[
\widehat{X}^{(k)}(x) = 
\sum_{\ell=1}^k \xi^{(k)}_\ell(x) \cdot x_\ell^{(k)}.
\]
Clearly we may assume that
$\xi^{(k)}_1,\dots,\xi^{(k)}_k$ are independent with respect to $\mu$.
Take independent copies $\Xh^{(k)}_1,\dots \Xh^{(k)}_n$ 
of $\Xh^{(k)}$ and
define the randomized algorithm $\Sh^{(k)}_n$ by \eqref{mceuler}.

\begin{theorem}\label{galg}
The algorithm $\Sh_N = \Sh_n^{(k)}$ with
$n = \lfloor N^{2/(2+\alpha)}\cdot 
(\ln N)^{-2(\alpha+\beta)/(2+\alpha)}\rfloor$ 
and 
$k=\lfloor N^{\alpha/(2+\alpha)}\cdot 
(\ln N)^{2(\alpha+\beta)/(2+\alpha)} \rfloor$ satisfies
\[
e(\Sh_N) \preceq N^{-1/(2+\alpha)}\cdot (\ln N)^{(\alpha+\beta)/(2+\alpha)}
\]
and
\[
\cost (\Sh_N) \preceq N.
\]
\end{theorem}

\begin{proof}
Proceed as in the proof of Theorem \ref{t8} to obtain
\[
e^2(\Sh_n^k) \preceq 1/n + (\ln k)^2 \cdot d_k^2 
\preceq 1/n +  k^{-2/\alpha}\cdot (\ln k)^{2(\alpha+\beta)/\alpha}
\]
by means of \eqref{n3} and Proposition \ref{pqwg}. The
asymptotic
upper bound for the error of $\Sh_N = \Sh_n^{(k)}$
now follows from the particular choice of $n$ and $k$.
Clearly, $\cost (\Sh_n^{(k)}) \asymp k \cdot n$.
\end{proof}

Combine Theorems \ref{tlg} and \ref{galg} to conclude that
the algorithm $\Sh_N$ is almost optimal in the following
sense.

\begin{cor}\label{c4}
There exists a constant $c >0$ such that
\[
e(\Sh_N) \leq c \cdot \eps^\ran_N \cdot
(\ln N)^{1+2/(\alpha(2+\alpha))}\cdot 
(\ln \ln N)^{-2\beta/(\alpha(2+\alpha))}
\]
holds for infinitely many integers $N$.
\end{cor}

\begin{rem}\label{r1000}
A slightly better upper bound is available
if the Banach space $\XX$ is
B-convex, e.g., if $\XX$ is an $L_p$-space with $p\in \left]1,\infty\right[$. 
Instead of \eqref{n3} we then have
\begin{equation}\label{n4}
\left( \int_{\XX} \| x - \widehat{X}^{(k)}(x) \|^2 \, 
\mu (dx) \right)^{1/2} \preceq d_k,
\end{equation}
see Creutzig (2002, Cor.\ 3.4.2), which yields
\[
e(\Sh_N) \preceq N^{-1/(2+\alpha)}\cdot (\ln N)^{\beta/(2+\alpha)}
\]
in Theorem \ref{galg}.
Both of the estimates \eqref{n3} and \eqref{n4} are proven
non-constuctively.

For a number of Gaussian measures on function spaces
the Karhunen-Lo\'eve expansion is explicitly known, and hereby we get
an approximation $\widehat{X}^{(k)}$ that satisfies \eqref{n4}, 
if $\XX$ is any $L_p$-space with $p\in \left]1,\infty\right[$.
In particular for an $L_2$-space $\XX$ and $\beta=0$ the
upper bound \eqref{n4}
is due to Wasilkowski, Wo\'zniakowski (1996, p.\ 2076).

Consider the distribution $\mu$ of
the $d$-dimensional fractional Brownian sheet with Hurst parameter
$H \in \left]0,1\right[$ on the space $\XX = C([0,1]^d)$.
In this case a direct approach yields
\begin{equation}\label{n5}
\left( \int_{\XX} \| x - \widehat{X}^{(k)}(x) \|^2 \, 
\mu (dx) \right)^{1/2} \preceq k^{-H} \cdot (\ln k)^{H(d-1) +d/2},
\end{equation}
see K\"uhn, Linde (2002).
See also
Ayache, Taqqu (2003) for a wavelet approximation
$\widehat{X}^{(k)}$ in the case $d=1$ and 
Dzhaparidze, van Zanten (2005) 
for a trigonometric approximation $\widehat{X}^{(k)}$
in the case $d \geq 1$, which both satisfy this estimate.
{}From \eqref{n5} we get
\[
e(\Sh_N) \preceq N^{-1/(2+1/H)}\cdot (\ln N)^{d/2 - 1/(2+1/H)}
\]
in Theorem \ref{galg}.
\end{rem}

\section*{Acknowledgments}

\noindent
We thank Jakob Creutzig and
Karsten Gro\ss e-Brauckmann for 
valuable discussions and comments.
In particular, Jakob 
pointed out the proof of Lemma \ref{ljakob} to us.

\section*{References}

{\small

\noindent
Ayache, A., Taqqu, M. S. (2003),
Rate optimality of wavelet series approximations of fractional Brownian
motion, \JFAA {\bf 9}, 451--471.
\medskip

\noindent
Bakhvalov, N. S. (1959),
On approximate computation of integrals (in Russian),
Vestnik MGV, Ser. Math. Mech. Astron. Phys. Chem. {\bf 4}, 3--18.
\medskip

\noindent
Bakhvalov, N. S. (1971),
On the optimality of linear methods for operator approximation
in convex classes of functions,
\CMMP {\bf 11}, 244--249.
\medskip

\noindent
Belinsky, E., Linde, W. (2002),
Small ball probabilities of fractional Brownian sheets via
fractional integration operators, \JTP {\bf 15}, 589--612.
\medskip

\noindent
Bouleau, N., L\'epingle, D. (1994),
Numerical Methods for Stochastic Processes,
Wiley, New York.%
\medskip\par%
\noindent
Cs\'aki, E. (1984), 
On small values of the square integral of a multiparameter 
Wiener process, in: Statistics and Probability, 
J.\ Mogyorodi, I.\ Vincze, W.\ Wertz,
eds., pp. 19--26, Reidel, Dordrecht.%
\medskip%
\par%
\noindent
Creutzig, J. (2002),
Approximation of Gaussian random vectors in Banach spaces,
Ph.D. Dissertation, Universit\"at Jena.
\medskip

\noindent
Dereich, S. (2003), High resolution coding of stochastic processes 
and small ball probabilities. Ph.D. Dissertation, TU Berlin.
\medskip

\noindent
Dereich, S. (2004),
The quantization complexity of diffusion processes, Preprint,
arXiv: \!math.PR/ 0411597.
\medskip

\noindent
Dereich, S., Fehringer, F., Matoussi, A., and Scheutzow, M.
(2003),
On the link between small ball probabilities and the quantization 
problem,
\JTP {\bf 16}, 249--265.
\medskip

\noindent
Dereich, S., Scheutzow, M. (2005),
High-resolution quantization and entropy coding for fractional 
Brownian motion, Preprint, arXiv: math.PR/0504480.
\medskip

\noindent
Dzhaparidze, K., van Zanten, H. (2005),
Optimality of an explicit series expansion of the fractional Brownian sheet,
\SPL {\bf 71}, 295--301.
\medskip

\noindent
Faure, O. (1992),
Simulation du mouvement brownien et des diffusions,
Th\`ese, ENPC Paris.
\medskip

\noindent
Fill, J. A., Torcaso, F.. (2004),
Asymptotic analysis via Mellin transforms for small deviations
in $L^2$-norm of integrated Brownian sheets,
\PTRF {\bf 130}, 259--288..
\medskip

\noindent
Graf, S., Luschgy, H. (2000),
Foundations of Quantization for Probability Distributions,
\LNM {\bf 1730}, Springer-Verlag, Berlin.
\medskip

\noindent
Gray, R.\ M., Neuhoff, D.\ L., Shields, P.\ C. (1975),
A generalization of Ornstein's $\overline{d}$ distance with
applications to information theory,
\AAP {\bf 3}, 315--328.
\medskip

\noindent
Kahane, J.-P. (1993),
Some Random Series of Functions,
Cambridge Univ.\ Press, Cambridge.
\medskip

\noindent
Kantorovich, L.\ V., Rubinstein, G.\ S. (1958),
On a space of completely additive functions (in Russian),
Vestnik Leningrad Univ. 13, no. 7, Ser. Mat. Astron. Phys.
{\bf 2}, 52--59. 
\medskip

\noindent
K\"uhn, T., Linde, W. (2002),
Optimal series representation of fractional Brownian sheets,
\BER {\bf 8}, 669--696.
\medskip

\noindent
Kuelbs, J., Li, W.\ V., Shao Q.\ M. (1995),
Small ball estimates for fractional Brownian motion under 
H\"older norm and Chung's functional LIL,
J. Theoret. Probab. {\bf 8}, 361--386.
\medskip

\noindent
Li, W.V.,  Shao, Q.-M. (1999), 
Small ball estimates for Gaussian processes under the Sobolev norm,
J. Theoret. Probab. {\bf12}, 699--720.
\medskip

\noindent
Li, W. V.,  Shao, Q.-M. (2001), 
Gaussian processes: inequalities, small ball probabilities and 
applications,
in: Stochastic Processes: Theory and Methods, 
Handbook of Statist., Vol.\ 19, D.\ N. Shanbhag, C.\ R.\ Rao, eds.,
pp.\ 533--597, North-Holland, Amsterdam.
\medskip

\noindent
Luschgy, H., Pag\`{e}s, G. (2003),
Functional quantization of 1-dimensional Brownian diffusion
processes, Preprint, Universit\'{e} de Paris VI, LPMA no. 853.
\medskip 

\noindent
Luschgy, H., Pag\`{e}s, G. (2004),
Sharp asymptotics of the functional quantization problem for 
Gaussian processes,
\AAP {\bf 32}, 1574--1599.
\medskip

\noindent
Maiorov, V. (1993),
Average $n$-widths of the Wiener space in the $L_\infty$-norm,
\JC {\bf 9}, 222--230.
\medskip

\noindent
Nemirovsky, A. S., Yudin, D. B. (1983),
Problem Complexity and Method Efficiency in Optimization,
Wiley, New York.
\medskip

\noindent
Novak, E. (1988),
Deterministic and Stochastic Error Bounds in Numerical Analysis,
\LNM {\bf 1349}, Springer-Verlag, Berlin.
\medskip

\noindent
Novak, E. (1995),
The real number model in numerical analysis,
\JC {\bf 11}, 57--73.%
\medskip\par%
\noindent
Pag\`es, G., Printems, J. (2004),
Functional quantization for pricing derivatives,
Preprint, Universit\'{e} de Paris VI, LPMA no. 930.
\medskip
 
\noindent
Ritter, K. (2000),
Average-Case Analysis of Numerical Problems,
\LNM {\bf 1733}, Springer-Verlag, Berlin.
\medskip

\noindent
Sakai, T. (1996),
Riemannian Geometry,
\TMM {\bf 149}, AMS, Rhode Island.
\medskip

\noindent
Smolyak, S. A. (1965),
On optimal restoration of functions and functionals of them
(in Russian), Candidate Dissertation, Moscow State University.
\medskip

\noindent
Traub, J. F., Wasilkowski, G. W., Wo\'zniakowski, H. (1988),
Information-Based Complexity,
Academic Press, New York.
\medskip

\noindent
Wasilkowski, G.\ W. (1989),
Randomization for continuous problems,
\JC {\bf 5}, 195--218.%
\medskip\par%
\noindent
Wasilkowski, G.\ W., Wo\'{z}niakowski, H. (1996),
On tractability of path integration,
\JMP {\bf 37}, 2071--2088.
\medskip

\noindent
Wasilkowski, G.\ W., Wo\'{z}niakowski, H. (2001),
Complexity of weighted approximation over $\R^d$,
\JC {\bf 17}, 722--740.
\medskip

}

\end{document}